\documentclass[11pt,a4paper]{article}
\usepackage[T1]{fontenc}
\usepackage[utf8]{inputenc}
\usepackage{lmodern}
\usepackage[margin=25mm]{geometry}
\usepackage{amsmath,amssymb}
\usepackage{microtype}
\usepackage{enumitem}
\usepackage{fancyvrb}
\usepackage{titlesec}
\usepackage[hidelinks]{hyperref}
\usepackage{xurl}
\hypersetup{pdftitle={A uniform double-factorial construction for consecutive-pattern packings},
 pdfauthor={Zolt\'an Nagy},pdfsubject={Consecutive permutation patterns and linear extensions}}
\titleformat{\section}{\normalfont\large\bfseries}{}{0pt}{}
\titlespacing*{\section}{0pt}{13pt plus 2pt minus 1pt}{5pt plus 1pt}
\RecustomVerbatimEnvironment{verbatim}{Verbatim}{fontsize=\small}
\title{\vspace{-1.5em}\Large\bfseries A uniform double-factorial construction\\for consecutive-pattern packings}
\author{Zolt\'an Nagy}
\date{26 October 2026}
\begin{document}
\maketitle
\vspace{-1.5em}
\textbf{Abstract.} For every overlap \(m\geq 3\) and every pattern
length \(K\geq 2m+1\), we construct an explicit permutation
\(p\in\mathfrak{S}_K\) such that the number of permutations containing
\(j\) prescribed occurrences of \(p\), with successive starts separated
by \(K-m\), is \((2j-1)!!\) for every \(j\geq 1\). Additional
occurrences are permitted. We prove that the corresponding occurrence
poset is the ordinal sum of a chain on \(m+j(K-m-2)\) elements and a
comb on \(2j\) elements. An insertion bijection counts the linear
extensions of the comb. This gives the uniform representatives requested
in the double-factorial part of Question 6 in Section 7 of Blitvić,
Kammoun and Steingrímsson, arXiv:2610.29846v1. A Lean 4 formalization
and its recorded verification data accompany the paper.


\section{1. Introduction}\label{introduction}

Blitvić, Kammoun and Steingrímsson {[}1{]} study permutations covered by
prescribed occurrences of a consecutive pattern with a fixed overlap.
They relate their enumeration to linear extensions of associated posets.
In Section 7, Question 6, they ask for a uniform family realizing the
sequence \(1,3,15,105,\ldots\), with one representative for every
overlap \(m>2\) and every length \(K>2m\).

We give such a family. Its overlapping boundary variables form
alternating rows that concatenate into a single chain. The maximum of
each prescribed occurrence supplies a leaf attached to the final row.
After removal of a forced initial chain, the resulting poset is a comb
whose number of linear extensions is an odd double factorial. The
construction varies with \(m\) and \(K\), but not with the number of
prescribed occurrences. The proof is self-contained; the formalization
is supplementary.

\section{2. Packings and the explicit
pattern}\label{packings-and-the-explicit-pattern}

Write \([n]=\{1,\ldots,n\}\) and let \(\mathfrak{S}_n\) denote the
permutations of \([n]\), written in one-line notation. Fix
\(1\leq m<K\), a pattern \(p\in\mathfrak{S}_K\) and an integer
\(j\geq 1\). Put
\[
d=K-m,\qquad n=K+(j-1)d=m+jd.
\]

Following {[}1, Definition 2.6 and Remark 2.7{]}, let \(N_{p,m}(j)\)
count the permutations \(s\in\mathfrak{S}_n\) for which
\[
p_h<p_k\quad\Longleftrightarrow\quad s_{rd+h}<s_{rd+k}
\qquad(0\leq r<j,\ 1\leq h,k\leq K).
\]

Thus the prescribed starts are \(1,1+d,\ldots,1+(j-1)d\). No
restrictions are imposed on occurrences at other starts, and each full
permutation \(s\) is counted once.

For the rest of the paper suppose that \(m\geq 3\) and \(K\geq 2m+1\).
Set \(c=K-2m\geq 1\), so that \(d=m+c\). For \(1\leq i\leq m\), let
\(\varepsilon_i=1\) when \(i\) is odd and \(\varepsilon_i=0\) when \(i\)
is even, and define
\[
a_i=c+2i-2+\varepsilon_i.
\]

Set \(b_1=1\) and, for \(2\leq i\leq m\), set
\(b_i=c+2i-1-\varepsilon_i\). Define
\[
p=p(m,K)=(a_1,\ldots,a_m,\ 2,3,\ldots,c,\ K,\ b_1,\ldots,b_m).
\]

The segment \(2,3,\ldots,c\) is empty when \(c=1\).

\textbf{Theorem 1.} \emph{For every \(m\geq 3\) and \(K\geq 2m+1\), the
word \(p(m,K)\) is a permutation in \(\mathfrak{S}_K\). For every
\(j\geq 1\) it satisfies}
\[
N_{p(m,K),m}(j)=(2j-1)!!=\prod_{t=1}^{j}(2t-1).
\]

\textbf{Lemma 2 (validity of the pattern).} \emph{The word \(p(m,K)\)
contains every element of \([K]\) exactly once.}

\emph{Proof.} Its length is \(m+(c-1)+1+m=K\). Since \(a_1=c+1\) and
\(b_1=1\), the first pair together with the padding uses precisely
\(1,\ldots,c+1\). For each \(i\geq 2\),
\[
\{a_i,b_i\}=\{c+2i-2,c+2i-1\}.
\]

These disjoint pairs use \(c+2,\ldots,K-1\). The separate entry \(K\)
completes \([K]\). \(\square\)

For example, the patterns for \((m,K)=(3,7),(3,8),(4,9)\) are,
respectively,
\[
(2,3,6,7,1,4,5),\qquad (3,4,7,2,8,1,5,6),
\]
\[
(2,3,6,7,9,1,4,5,8).
\]

\section{3. From prescribed occurrences to a
poset}\label{from-prescribed-occurrences-to-a-poset}

Introduce a vertex \(v_t\) for each position \(t\in[n]\). For every
prescribed block \(r\), impose the chain
\[
v_{rd+p^{-1}(1)}<v_{rd+p^{-1}(2)}<\cdots<v_{rd+p^{-1}(K)}.
\]

We use consecutive comparisons in these chains as generating edges.
Their transitive closure contains all the comparisons within each
prescribed block. When this directed graph is acyclic, its closure
defines a strict partial order, denoted \(P_{p,m}(j)\). Acyclicity for
the present construction will follow from Proposition 4, rather than
being assumed. For a finite poset \(P\), write \(e(P)\) for its number
of linear extensions.

\textbf{Lemma 3 (the counting correspondence).} \emph{If the occurrence
graph is acyclic, then \(N_{p,m}(j)=e(P_{p,m}(j))\).}

\emph{Proof.} Given a linear extension \(v_{t_1}<\cdots<v_{t_n}\),
assign \(s_{t_i}=i\). This is a permutation of \([n]\) and satisfies all
the prescribed block comparisons. Conversely, ordering the positions of
an admissible permutation by their values gives a linear extension.
These operations are inverse. Since each block is a total order on
distinct values, its generating chain is equivalent to the full
relative-order condition in Section 2. \(\square\)

\section{4. The exact chain-and-comb
decomposition}\label{the-exact-chain-and-comb-decomposition}

Let \(C_h\) be a chain on \(h\) elements. For posets \(A\) and \(B\),
their ordinal sum \(A\oplus B\) retains the orders within the two parts
and declares every element of \(A\) smaller than every element of \(B\).

Define \(Q_j\) on distinct elements
\(q_1,\ldots,q_j,\ell_1,\ldots,\ell_j\) by the generating relations
\[
q_1<q_2<\cdots<q_j,\qquad q_t<\ell_t\quad(1\leq t\leq j).
\]

There are no other generating relations. In particular, \(q_s<\ell_t\)
exactly when \(s\leq t\), and the leaves \(\ell_t\) are pairwise
incomparable. Let \(Q_0\) be the empty poset.

\textbf{Proposition 4.} \emph{For \(p=p(m,K)\) and every \(j\geq 1\),
the occurrence graph is acyclic and}
\[
P_{p,m}(j)\cong C_{\,m+j(d-2)}\oplus Q_j,
\qquad d=K-m.
\]

\emph{Proof.} Name the position vertices as follows:
\[
x_{r,i}=v_{rd+i}\qquad(0\leq r\leq j,\ 1\leq i\leq m),
\]
\[
u_{r,t}=v_{rd+m+t}\qquad(0\leq r<j,\ 1\leq t<c),
\qquad y_r=v_{(r+1)d}\qquad(0\leq r<j).
\]

For each \(r<j\), the vertices \(x_{r,1},\ldots,x_{r,m}\), the
\(u_{r,t}\) and \(y_r\) occupy exactly the interval of positions
\(rd+1,\ldots,(r+1)d\). The remaining vertices
\(x_{j,1},\ldots,x_{j,m}\) occupy \(jd+1,\ldots,jd+m\). Hence these
names are distinct and exhaust all \(n=m+jd\) vertices.

Block \(r\) consists, in position order, of its \(x_{r,i}\), its
\(u_{r,t}\), \(y_r\), and the \(x_{r+1,i}\). Since \(d>m\), nonadjacent
prescribed blocks do not intersect, and each \(u_{r,t}\) and \(y_r\)
belongs to only one block.

Define the lower and upper vertex in row \(i\) of block \(r\) by
\[
(L_{r,i},H_{r,i})=(x_{r+1,i},x_{r,i})\quad\text{for odd }i,
\]
\[
(L_{r,i},H_{r,i})=(x_{r,i},x_{r+1,i})\quad\text{for even }i.
\]

The ranks assigned by \(p\) to \(L_{r,1}\), to the padding vertices, and
to \(H_{r,1}\) are \(1\), then \(2,\ldots,c\), then \(c+1\). For
\(i\geq 2\), the ranks at \(L_{r,i},H_{r,i}\) are \(c+2i-2,c+2i-1\),
respectively. Finally, \(y_r\) has rank \(K\). Thus the defining chain
of block \(r\) starts with
\[
L_{r,1}<u_{r,1}<\cdots<u_{r,c-1}<H_{r,1},
\]

and continues with
\[
H_{r,1}<L_{r,2}<H_{r,2}<\cdots<L_{r,m}<H_{r,m}<y_r.
\]

Together these two segments give the complete sorted block order. The
padding segment is omitted when \(c=1\).

Within an odd row \(i\) these chains give
\(x_{j,i}<x_{j-1,i}<\cdots<x_{0,i}\); within an even row they give
\(x_{0,i}<x_{1,i}<\cdots<x_{j,i}\). The rows are completely separated,
not merely partially interlaced. Indeed, for odd \(i<m\), block \(0\)
supplies
\[
\max(\text{row }i)=x_{0,i}<x_{0,i+1}=\min(\text{row }i+1).
\]

For even \(i<m\), block \(j-1\) supplies
\[
\max(\text{row }i)=x_{j,i}<x_{j,i+1}=\min(\text{row }i+1).
\]

The padding vertices subdivide the intervals from \(x_{r+1,1}\) to
\(x_{r,1}\). Therefore all the \(x\)- and \(u\)-vertices form a chain,
denoted \(B\).

To establish the exact order, define \(R\) to be the poset consisting of
the chain \(B\) together with the generating relations \(H_{r,m}<y_r\),
and no outgoing relations from any \(y_r\). This is acyclic. Every
generating relation of \(R\) follows from the block chains: the
within-row relations and padding relations occur there, the row joins
were just proved, and each leaf edge is the last edge of its block.
Conversely, every comparison in every block chain holds in \(R\),
because the rows of \(B\) are in order and the required final leaf edge
is present. The two transitive closures therefore coincide. In
particular, the occurrence graph is acyclic, and there are no extra
comparabilities. The predecessors of \(y_r\) are exactly the elements of
\(B\) at or below \(H_{r,m}\).

Write the final row of \(B\) in increasing order as
\(q_0<q_1<\cdots<q_j\). For odd \(m\), use
\[
q_t=x_{j-t,m}\ (0\leq t\leq j),\qquad
\ell_t=y_{j-t}\ (1\leq t\leq j).
\]

For even \(m\), use
\[
q_t=x_{t,m}\ (0\leq t\leq j),\qquad
\ell_t=y_{t-1}\ (1\leq t\leq j).
\]

In either parity, the leaf edges become precisely \(q_t<\ell_t\). All
elements of \(B\) through \(q_0\) precede every remaining vertex.
Removing this initial chain leaves exactly \(Q_j\). There are
\(n-2j=m+j(d-2)\) vertices in the initial chain, proving the asserted
ordinal-sum decomposition. \(\square\)

\section{5. Enumeration of the comb and proof of Theorem
1}\label{enumeration-of-the-comb-and-proof-of-theorem-1}

\textbf{Lemma 5.} \emph{For every \(j\geq 0\),}
\[
e(Q_j)=\prod_{t=1}^{j}(2t-1),
\]

\emph{where the empty product equals \(1\).}

\emph{Proof.} The empty poset has one linear extension. Suppose
\(j\geq 1\). In every extension of \(Q_j\), the element \(q_1\) comes
first: it lies below the whole spine and hence below every leaf. Delete
\(q_1\) and \(\ell_1\), and decrease the remaining subscripts by one.
The resulting word is an extension of \(Q_{j-1}\).

Conversely, start with an extension of \(Q_{j-1}\), increase its
subscripts by one, insert \(\ell_1\) into any of the \(2j-1\) gaps of
this word of length \(2j-2\), and prepend \(q_1\). Each insertion is
valid because \(\ell_1\) is incomparable with every remaining element
once \(q_1\) has been deleted. Deletion recovers both the smaller
extension and the chosen gap. Consequently,
\[
e(Q_j)=(2j-1)e(Q_{j-1}),
\]

which proves the product formula by induction. \(\square\)

\emph{Proof of Theorem 1.} Lemma 2 proves that \(p(m,K)\) is a
permutation. Proposition 4 and Lemma 3 identify the required count with
the number of extensions of \(C_{m+j(d-2)}\oplus Q_j\). Every extension
begins with the forced chain in its unique order, so its removal gives a
bijection with the extensions of \(Q_j\). Lemma 5 now gives the result.
The formula defining \(p\) uses only \(m\) and \(K\), so the same
pattern works for all \(j\geq 1\). \(\square\)

There is also a direct bijection between extensions of \(Q_j\) and
perfect matchings of \([2j]\). Pair the positions occupied by \(q_t\)
and \(\ell_t\). Conversely, order the pairs of a matching by their
smaller endpoints, label them \(1,\ldots,j\), and place \(q_t\) at the
smaller endpoint and \(\ell_t\) at the larger. These maps are inverse.
Choosing the partner of the smallest remaining position recursively
gives the same odd double factorial.

The boundary \(K=2m+1\) requires no separate argument: the padding is
empty and all row joins remain valid. When \(j=1\), the prefix has
\(K-2\) vertices and \(Q_1\) is a two-element chain, giving one packing.
The construction settles the requested double-factorial representatives;
it does not classify all equinumerous patterns or address the other
enumeration sequences discussed in {[}1, Section 7{]}.

\section{6. Formalization and
reproducibility}\label{formalization-and-reproducibility}

The accompanying archive \texttt{PatternPacking-source.tar.gz} contains
a Lean 4/Mathlib formalization. The original counting predicate is
defined directly on permutations by the blockwise relative-order
biconditionals; it is not defined using the comb or the expected count.
The main declaration, written here with ASCII alternatives for Lean's
natural-number and inequality notation, is

\begin{verbatim}
theorem main_theorem (m K j : Nat)
    (hm : 3 <= m) (hK : 2*m+1 <= K) (hj : 1 <= j) :
    Nat.card (Packings (candidate m K hm hK) m j)
      = Comb.oddDF j
\end{verbatim}

It appears in \texttt{PatternPacking/Main.lean}. The declaration
\texttt{occurrence\_structure} in \texttt{StructuralIsomorphism.lean}
expresses the exact structural correspondence. The comb bijection and
its enumeration are in \texttt{Comb.lean}; the permutation-to-extension
bridge is in \texttt{Bridge.lean}.

The archive pins Lean to version 4.35.0 and Mathlib to revision

\begin{verbatim}
56a005855ea17e0b0958f73ec801ae140b0a1679
\end{verbatim}

Its supplied verification records report successful project and
fresh-source builds, with shared dependency caches, and list only
\texttt{propext}, \texttt{Classical.choice} and \texttt{Quot.sound} as
axioms of the main theorem. These are standard foundational axioms of
Lean {[}2{]}. The records also report 290 supplementary regression
cases. These finite tests are not premises of the proof.

These are the recorded checks shipped with the formalization; no
separately reproduced build or standalone kernel replay is claimed here.
After installing the pinned toolchain, the core reproduction commands,
run at the extracted project root, are

\begin{verbatim}
lake exe cache get
lake build
lake env lean PatternPacking/Main.lean
lake env lean PatternPacking/Audit.lean
\end{verbatim}

The archive includes its dependency manifest, source hashes and
verification logs. The mathematical proof in Sections 2--5 does not
depend on these execution records. Formal verification of a specified
statement and establishing its originality are distinct questions.

\section{AI-use disclosure}\label{ai-use-disclosure}

ChatGPT, Codex and Claude were used substantively in developing the
construction and proof, implementing the Lean formalization, adversarial
checking, and preparing this manuscript. These tools are not authors.
Responsibility for the submitted work rests with the human author.

\section*{References}
\begin{enumerate}[label={[\arabic*]},leftmargin=1.7em,itemsep=5pt]
\item Natasha Blitvi\'c, Mohammed Slim Kammoun and Einar Steingr\'imsson.
\emph{Permutations, Pattern Covers, Posets and Positivity: Part I: Patterns and Posets}.
arXiv:2610.29846v1, 2026. \url{https://arxiv.org/abs/2610.29846v1}.
\item The Lean development team. \emph{The Lean Language Reference},
``Validating a Lean Proof.''
\url{https://lean-lang.org/doc/reference/latest/ValidatingProofs/}.
Accessed 26 October 2026.
\end{enumerate}

\end{document}